\RequirePackage{fix-cm}
\documentclass[draft]{svjour3}                     
\smartqed  
\usepackage{graphicx}

\usepackage{amsmath}
\usepackage{mathtools}
\usepackage{amssymb}

\newcommand{\be}{\begin{equation}}
\newcommand{\ee}{\end{equation}}
\begin{document}

\title{Equivalence of variational principles to determine the speed of scalar reaction diffusion fronts}

\titlerunning{On the speed of scalar reaction diffusion fronts }        

\author{R. D. Benguria       \and
        M. C. Depassier 
}


\institute{R. D. Benguria \& M. C. Depassier \at
              Instituto de F\'\i sica\\
              Pontificia Universidad Cat\'olica de Chile\\
              Casilla 306, Santiago 22, Chile \\
              \email{rbenguri@uc.cl,   mcdepass@uc.cl}           
}

\date{Received: date / Accepted: date}

\maketitle

\begin{abstract}
The determination of the  speed of travelling fronts of the scalar reaction diffusion equation has been the subject of much study.  Using different approaches seemingly disconnected variational principles have been established. 
 The purpose of this work is to show the connection between them. 
For monostable reaction terms, we prove that a principle established by  Hadeler and Rothe in 1975 and a second one  by Benguria and Depassier in 1996 are logically equivalent, that is, either can be derived from the other.  Two variational principles,  formulated for arbitrary reaction terms,  are shown to be related by a suitable change of variables. Finally a variational principle proven  for monostable reaction terms is shown to be  a formulation of the two 
previous  ones in yet another independent variable.

\end{abstract}

\section{Introduction}
\label{intro}
The time evolution of the scalar reaction diffusion equation with constant diffusivity
\be\label{RD1}
u_t =  u_{xx}+ f(u) 
\ee
where the nonlinear reaction term has at least two equilibria, which we may  set  at $u=0$ and $u=1$, that is  $f(0) = f(1) = 0$ has been the subject of much study since its introduction  by Fisher  \cite{Fisher1937} in the study of the propagation of an advantageous  gene through a population,  and the proof by  Kolmogorov, Petrovsky and Piskunov  \cite{Kolmogorov1937}  (KPP) that the speed of the traveling wave of the now called Fisher-KPP equation
\be
u_t = u_{xx} + u(1-u)
\ee
is $c=2.$  More generally,  KPP  showed that a sufficient localized initial condition $u(x,0)$ evolves into a monotonic traveling wave $u(x,t) = q(x - c  t)$  joining the stable state $u=1$ to the unstable state $u=0$. The speed is the minimal speed for which a monotonic front exists  and it is given by  $c\equiv c_{\rm{KPP}} = 2 \sqrt{f'(0)} $ for any monostable reaction term, that is, for any  $f(u)$ such that $f( 0) = f(1) =0, f(u)>0 \in (0,1),  f'(0) >0$ which in addition satisfies $f(u) <  f'(0) u.$  While the equation proposed by Fisher satisfies these hypothesis and the speed of the fronts is $c=2$, reaction terms which appear in the context of combustion theory do not, and Zeldovich and Frank-Kameneetskii \cite{Zeldovich1938} studied a reaction term concentrated at $u=1$ and vanishing elsewhere and showed that in this case  $c = c_{\rm{ZFK}}$ with
\be 
c_{\rm{ZFK}} = \sqrt{ 2 \int_0^1 f(u) \, du}.
\ee 
This asymptotic value of the speed was later shown \cite{Berestycki1992}  to be a lower bound for the speed for any monostable reaction term. 

 The time evolution for general reaction terms was studied by Aronson and Weinberger \cite{Aronson1975,Aronson1978} who showed that sufficiently localized initial conditions evolve into a front  which propagates with speed $c_0$ such that 
 \be\label{ulb}
  2 \sqrt{f'(0)} \le  c_0 \le 2 \sqrt{\sup_{0 \le u\le 1}(f(u)/u)}.
  \ee
   The asymptotic speed of propagation is the minimal speed for which a monotonic front joining the stable to unstable equilibrium point exists.

A monotonic decaying traveling front $u(x - c t)$ obeys the differential equation
$$
u_{zz} + c  u_z + f(u) = 0,
$$
with
$$
 \lim_{z \rightarrow -\infty }u(z) = 1, \quad 
\lim_{z \rightarrow \infty }u(z) = 0,  \quad u_z < 0,
$$
where $z = x -c t$. 
It is convenient to work in phase space;  defining as usual $p(u) = - q_z$, the problem reduces to finding the solutions of minimal speed of 
\begin{equation}
p(u)\, \frac{d p(u) }{du} - c  p(u) + f(u) = 0 , \qquad \qquad p(0) = p(1) = 0 , \qquad {\rm and}\,\, p(u)> 0.
\label{phase}
\end{equation}

When the upper and lower bounds in (\ref{ulb}) do not coincide the speed can be determined either for particular exactly solvable cases or estimated with any desired accuracy using variational principles.  The purpose of this work is to 
show the connection between these variational principles.

Three variational principles have been established for monostable reaction terms satisfying   $f(u) >0 \in (0,1)$, $f(0) = f(1) =0$ and $f'(0) >0$.

The first principle (VP1) , formulated by Hadeler and Rothe in 1975 \cite{Hadeler1975}, states that   the minimal speed  satisfies 
\be \label{vp1}
c_0 = \inf_{\alpha} \sup_{0\le u \le 1} \left( \alpha'(u) + \frac{f(u)}{\alpha(u)} \right),  
\ee
where $\alpha \in C^1([0,1])$,  $\alpha(0) = 0, \alpha'(0) > 0$ and  $\alpha(u) > 0$ in $(0,1)$. From this variational principle, upper bounds on the speed are obtained with suitable trial functions $\alpha(u)$, that is, 
$$
c_0 \le \sup_{0\le u \le 1} \left( \alpha'(u) + \frac{f(u)}{\alpha(u)} \right). 
$$

The second principle (VP2), formulated by Benguria and Depassier in 1996 \cite{BeDeCMP1996}, states that 
\be\label{vp2}
c_0 = \sup_g \frac{ 2 \int_0^1 \sqrt{ f(u) g(u)h(u)} \, du}{\int_0^1 
 g(u) \, du },
\ee
where $g(u) \in C^1((0,1])$ is any positive decreasing function for which the integrals exist and $h(u) = - g'(u) >0$.  The supremum is 
attained when $c_0 > c_{\rm{KPP}}$ for a special trial function $\hat g$ which satisfies $\hat g (1) =0$. From this variational principle lower bounds on the speed  are obtained,
$$
c_0 \ge  \frac{ 2 \int_0^1 \sqrt{ f(u) g(u)h(u)}\, du}{\int_0^1 
 g(u) \, du }.
$$
A third principle  (VP3) proven only  for  monostable reaction terms,  formulated by Lucia, Muratov and Novaga in 2004 \cite{Lucia2004}, states that $c_0$ is greater than $c_{\rm{KPP}}$ if there exists $c> c_{\rm{KPP}
}$ and a function $u \in  H_c^1(\mathbb{R})$ such that 
\be \label{vp3}
\Phi_c[u] = \int_{-\infty}^{\infty} \text{e}^{c z} \left(\frac{1}{2} u_z^2 - \int_0^u f(s) ds \right) \,dz  \le 0.
\ee

The reaction diffusion equation is also of interest in problems where the reaction term is not monostable. Monotonic traveling fronts arise for bistable reaction terms which satisfy $f(0)=f(1)=0$, and
 \be
 \label{caseB}
   f(u)  < 0 \text{ for } u \text{ on } (0, a),\,   f(u) >0 \text{ on } (a, 1),  \text{ and } \int_0^1 f(u) \, du >0.
 \ee
Monotonic traveling fronts also arise in combustion problems where  the reaction term satisfies $f(0)=f(1)=0$, and
 \be\label{caseC}
 f=0 \text{ on }  (0,a), \text{ and }  f(u) > 0 \text{ on }  (a, 1).
 \ee
 Another case of interest, encountered  when density dependent diffusion is present,   leads to the case 
  \be \label{caseD}
  f'(0) = 0, \, f(u) >0  \text{ on }  (0,1).
  \ee

A variational principle (VP4) formulated in phase space valid for any  reaction term $f(u)$ for which (\ref{RD1}) evolves into monotonic fronts, established in (\cite{BeDePRL1996}), states that
\be\label{vp4}
c_0^2 = \sup_g \frac{ 2 \int_0^1 f(q) g(q) dq}{  \int_0^1 g^2(q)/h(q)\,  dq},
\ee
where the supremum is taken over all positive decaying functions $g(u) \in C^1((0,1])$  for which the integrals exist. As before, $h(u) = - g'(u) >0$.    Whenever $c_0 \neq c_{\rm{KPP}}$ the supremum is attained for a trial function $\hat g$   given by 
\be\label{hatg}
\frac{\hat g'}{\hat g} = - \frac{c}{p}
\ee
and yields the minimal speed $c_0$ when  $c$ and $p(u)$ are the solution of (\ref{phase}). In this case $\hat g(1)=0$. 

Earlier and motivated by combustion problems, Rosen \cite{Rosen1958} in 1958, using standard calculus of variations, formulated a different principle, with no restrictions on the reaction term.  This principle (VP5), states that the wave profile $u(\xi)$  must be such that the action functional
\be\label{vp5}
\Lambda[u] = \frac{\int_0^\infty  \frac{1}{2}  u'(\xi)^2 d\xi}{\int_0^{\infty} \frac{1}{\xi^2} V d\xi}
\ee
assumes its minimum. 
Here, 
\be\label{V}
V(u)  =  \int_0^u f(s) ds. 
\ee
 This minimum is the  square of the inverse of the speed. Since the focus was on combustion problems no discussion on the validity of this principle for the KPP case was addressed. 

The purpose of this work is to show the connection between these variational principles. In Section 2 we show that VP1 $\Leftrightarrow$  VP2. In section 3 we show that the formulation on VP4 in coordinate space implies VP3. We show by a suitable coordinate change that VP4 and VP5  are identical.  While both VP4 and VP5 hold for arbitrary reaction terms,  VP5 did  not include the KPP case as it was originally stated.   VP4  has been extended to include monostable discontinuous reaction terms \cite{BeDeLo2012}.

\section{The relation between VP1 and VP2}

In this section we prove that each of these two variational principles implies the other. 

\subsection{VP2 implies VP1}

We assume that  (\ref{vp2}) holds and we will make use of   the inequality $2 a b \le \alpha a^2 + b^2/\alpha$ for any  $\alpha>0$, $\alpha \in C_1([0,1]).$ Choosing $a= \sqrt{h(u)}, b = \sqrt{f(u)g(u)}$ and $\alpha(u) $ such that $\alpha(0)=0, \alpha'(0) >0$ and $\alpha(u) >0 $ in $(0,1)$  we obtain from  (\ref{vp2})
 \begin{align*}
\begin{split}
c_0 =& \sup_g \frac{ 2 \int_0^1 \sqrt{ f(u) g(u)h(u)} \, du}{\int_0^1 
 g(u) \, du } \le  \sup_g \frac{  \int_0^1\left(  \alpha(u) h(u) + \frac{ f(u) g(u)}{\alpha(u)}  \, du \right) }{\int_0^1 
 g(u) \, du } \\
  \le&  \sup_g \frac{  \int_0^1\left(  \alpha'(u)  + \frac{ f(u) }{\alpha(u) }\right) g(u)  \,  du  }{\int_0^1 
 g(u) \, du } \le  \sup_u \left(  \alpha'(u)  + \frac{ f(u) }{\alpha(u) }\right), 
 \end{split}
\end{align*}
 where we integrated by parts and used the properties of $\alpha(u)$. This upper bound holds for any $\alpha(u)$ so we choose the function $\alpha(u)$ that gives the lowest upper bound, that is,
 $$
 c_0 \le \inf_{\alpha} \sup_u \left( \alpha'(u) + \frac{f(u)}{\alpha(u)} \right).
 $$
 Moreover it follows from (\ref{phase})  that when $\alpha(u) = p(u)$ one obtains the exact value of the speed, therefore  the first variational principle (\ref{vp1})  holds.

\subsection{VP1 implies VP2}

We start now from (\ref{vp1}).
 Since for any positive $g(u)$ we know that
$$
\frac{\int_0^1 \phi(u) g(u) \, du}{\int_0^1 g(u) \, du} \le \sup_u \phi(u),
$$
holds, 
it follows from 
(\ref{vp1})  that 
$$
c_0 = \inf_{\alpha} \sup_u \left( \alpha'(u) + \frac{f(u)}{\alpha(u)} \right) \ge  \inf_{\alpha}\frac{ \int_0^1 ( \alpha'(u) + f(u)\alpha(u)) g(u)  \, du}{\int_0^1 g(u) \,  du}
$$

Choosing the function $g(u)$ such that $g(1) = 0$ and $h(u)  = -g'(u)$, integrating by parts. we obtain
$$
c_0  \ge \inf_{\alpha} 
\frac{\int_0^1  \left( \alpha(u) h(u)  + \frac{f(u) g(u)  }
{\alpha(u)} \right) \, du}{\int_0^1 g(u) \, du} 
\ge  2 
\frac{ \int_0^1\sqrt{f(u) g(u) h(u)} \, du}{\int_0^1 g(u) \, du}.
$$
As shown in \cite{BeDeCMP1996} it can be proved that the supremum over $g$  of the right side above is precisely $c_0$.

\section{Reformulation of VP4 and its relation  to VP3 and VP5.}

As stated above, when the speed is not $c_{\rm{KPP}}$ the supremum in (\ref{vp4}) is attained for a trial function $\hat g(u)$ given by
$ \hat g'/ \hat g = - c/p$.  This value of $\hat g$ arises from maximizing (\cite{BeDePRL1996}) at fixed $u$ the function $\varphi(p) = c p(u) g(u) - p^2 (u) h(u)/2$. This inequality  holds for any $c>0$, $p(u)>0$ and the maximum in the variational principle is attained when $c= c_0$ and $p(u)$ is the corresponding solution of (\ref{phase}). 
Therefore we know that when $c_0 > c_{\rm{KPP}}$
$$
c_0^2 =  \frac{ 2 \int_0^1 f( q) \hat g(q) dq}{  \int_0^1 \hat g^2(q)/\hat h(q)\,  dq}.
$$

Next we formulate this principle in space variables.  Since $p(u) = - u_z$  (\ref{hatg}) can be integrated and yields
$$
\hat g (z) = \text{e}^{c z}.
$$
Recalling that $u= 0 $ at $z = +\infty$ and $u=1$ at $z = -\infty$, a straightforward change of variables yields
\be\label{down}
\int_0^1 \frac{\hat g^2(q)}{\hat h(q) } dq =  \frac{1}{c} \int_{-\infty}^{+\infty} \text{e}^{c z} (u_z)^2 dz.
\ee
  The numerator in (\ref{vp4}) can be written as 
\begin{align}
\begin{split}\label{up}
\int_0^1 f(q) \hat g(q) dq =& -\int_0^1 V(q)
 \frac{d \hat g}{d q} dq  
 = \int_{-\infty}^{+\infty} V(u(z))   \frac{d \hat g}{d z} dz \\
=&   c  \int_{-\infty}^{+\infty} \text{e}^{c z}  V(u(z))  dz. 
\end{split}
\end{align}

where $V$ is defined in (\ref{V}) and we used that $ V(0) =0, \hat g(1) =0$. 

 Replacing  (\ref{down}) and (\ref{up}) in (\ref{vp4}), we have that
\be \label{phasespace}
c_0^2 =  c^2  X_c 
 \ee
where we have defined
$$
X_c = \left(  \int_{-\infty}^{\infty} \text{e}^{c z} ( V(u(z))  dz\right) \Bigg/ \left( \int_{-\infty}^{\infty} \text{e}^{c z}  \frac{u_z^2}{2} dz\right).
$$
Therefore,  $X_c\ge 1$ for any $c$ such that $c_{\rm{KPP}} < c \le c_0$  or equivalently
\be 
\Phi_c = \int_{-\infty}^{\infty} \text{e}^{c z} \left(\frac{1}{2} u_z^2 -  V(u(z)) \right) dz  \le 0.
\ee
 Notice that $X_c=1$   when the trial function coincides with the minimal speed solution. The variational principle (\ref{vp5}) holds for any type of reaction term, not only for monostable reaction terms, and the optimizing $g$, $\hat g,$ exists for any reaction terms except in the KPP case  $c_0 = c_{\rm{KPP}}$.  In all cases at the solution $X_c=1$, and $\Phi_c=0$.   For monostable terms with $f'(0) = 0$, we know then that $X_c\ge 1$ ( $\Phi_c \le 0$) for any $0 < c< c_0$. 
 
 In order to establish the connection between VP4 with VP5 we use its reformulation \cite{BeDePRE2007} using as an independent variable $s = 1/g(u)$. As stated above, for $C^1([0,1])$ reaction terms the supremum in VP4 is attained for a trial function $\hat g$ which satisfies $\hat g(1) =0$.  Therefore without loss of generality one may restrict the trial function to satisfy $g(1) =0$.   Using this new independent variable,  (\ref{vp4}) can be written as
  \be\label{VP4-2}
 c_0^2 = \sup_{u(s)} \frac{\int_0^{\infty} V(u(s))/s^2}{\int_0^{\infty} (du/ds)^2 ds}
\ee
 where $V$ is given by (\ref{V}) and  the supremum is taken over positive increasing functions such that $u(0)=0$,  $\lim_{s\rightarrow\infty} u(s) =1$ and all integrals exist.   As proven in \cite{BeDeLo2012} the maximizer in ( \ref{VP4-2})  is not attained for reaction terms which satisfy the KPP criterion $f(u) \le f'(0) u$.   The extension of this variational principle to positive discontinuous reaction terms,   such as  monostable reaction terms with a cut-off, has been proven in \cite{BeDeLo2012} . This form of VP4 coincides with the principle formulated by Rosen except for the KPP case which is not included when (\ref{vp5}) is formulated considering the maximum instead of the supremum.

 \section{Summary}
 
We have presented a unified view of five seemingly disconnected variational principles that have been formulated to characterize the speed of a propagating front of the scalar  reaction diffusion equation. 
 For monostable reaction terms we have shown the logical equivalence of the two  principles  VP1 and VP2 \cite{Hadeler1975,BeDeCMP1996}.  The variational principles VP4 and VP5 which are  valid for general reaction terms  are equivalent and VP3 can be derived from them.

 \begin{acknowledgements}
 This work was partially supported by Fondecyt (Chile) project \# 120-1055.

\end{acknowledgements}


 

 \end{document}